\documentclass[draft]{amsart}

\usepackage{color}

\usepackage{hyperref}
\hypersetup{colorlinks=true, linkcolor = blue}

\usepackage{pgf,tikz,graphicx}
\usepackage{amssymb}
\usepackage{enumerate}
\usepackage{stmaryrd}

\renewcommand{\P}{\mathcal{P}}
\newcommand{\D}{\mathcal{D}}
\renewcommand{\O}{\mathcal{O}}
\newcommand{\E}{\mathcal{E}}

\numberwithin{equation}{section}
\theoremstyle{plain}
\newtheorem{theorem}{Theorem}[section]

\newtheorem{corollary}[theorem]{Corollary}

\begin{document}

\makeatletter
\def\imod#1{\allowbreak\mkern10mu({\operator@font mod}\,\,#1)}
\makeatother

\author{Alexander Berkovich}
   \address{[A.B.] Department of Mathematics, University of Florida, 358 Little Hall, Gainesville FL 32611, USA}
   \email{alexb@ufl.edu}

\author{Ali Kemal Uncu}
\address{[A.K.U.] University of Bath, Faculty of Science, Department of Computer Science, Bath, BA2\,7AY, UK}
\email{aku21@bath.ac.uk}

\address{[A.K.U.] Johann Radon Institute for Computational and Applied Mathematics, Austrian Academy of Science, Altenbergerstraße 69, A-4040 Linz, Austria}
\email{akuncu@ricam.oeaw.ac.at}

\title[\scalebox{.9}{On Finite Analogs of Schmidt's Problem and Its Variants}]{On Finite Analogs of Schmidt's Problem and Its Variants}
     
\begin{abstract} 
We refine Schmidt's problem and a partition identity related to 2-color partitions which we will refer to as Uncu--Andrews--Paule theorem. We will approach the problem using Boulet--Stanley weights and a formula on Rogers--Szeg\H{o} polynomials by Berkovich--Warnaar, and present various Schmidt's problem alike theorems and their refinements. Our new Schmidt type results include the use of even-indexed parts' sums, alternating sum of parts, and hook lengths as well as the odd-indexed parts' sum which appears in the original Schmidt's problem. We also translate some of our Schmidt's problem alike relations to weighted partition counts with multiplicative weights in relation to Rogers--Ramanujan partitions.
\end{abstract}   
   
 
 \thanks{Research of the second author is partly supported by EPSRC grant number EP/T015713/1 and partly by FWF grant P-34501N}
 

\date{\today}
   
\maketitle

\section{Introduction}\label{section1}
%%%%%%%%%%%%%%%%%%%%%%%%%%%%%%%%%%%%%%%%%%%%%%%%%%%%%%%%%%%%%%%%%%%%%%%%%%%%%%%%%%%%%%%%%%%%%%%%%%%%%

A \textit{partition} $\pi$ is a non-increasing finite sequence $\pi = (\lambda_1,\lambda_2,\dots)$ of positive integers. The elements $\lambda_i$ that appear in the sequence $\pi$ are called \textit{parts} of $\pi$. The \textit{total number of parts} is denoted by $\#(\pi)$. For positive integers $i$, we call $\lambda_{2i-1}$ odd-indexed parts, and $\lambda_{2i}$ even indexed parts of $\pi$. The sum of all the parts of a partition $\pi$ is called the \textit{size} of this partition and denoted by $|\pi|$. We say $\pi$ is a partition of $n$, if its size is $n$. The empty sequence $\emptyset$ is considered as the unique partition of zero.

Let $\D$ be the set of partitions into distinct parts, and let $\P$ be the set of all the (ordinary) partitions. The bivariate generating functions for these sets, where the exponent of $x$ is counting the number of parts and the exponent of $q$ is keeping track of the size of the partitions are known to have infinite product representations. Explicitly
\[\sum_{\pi\in\D} x^{\#(\pi)} q^{|\pi|} = (-xq;q)_\infty\quad\text{and}\quad\sum_{\pi\in\P} x^{\#(\pi)} q^{|\pi|} = \frac{1}{(xq;q)_\infty},\] where \[(a;q)_l := \prod_{i=0}^{l-1}(1-a q^i),\quad \text{for}\ \ l\in\mathbb{Z}_{\geq0}\cup \{\infty\}\] the standard $q$-Pochhammer symbol \cite{theoryofpartitions}.

In 1999, although the original problem was submitted in 1997, F. Schmidt \cite{Schmidt} shared his curious observations on a connection between partitions into distinct parts and ordinary partitions. For a given partition $\pi = (\lambda_1,\lambda_2,\dots)$, let \[\O(\pi):= \lambda_1+\lambda_3+\lambda_5+\dots,\] be the \textit{sum of odd-indexed parts}. Then, Schmidt observed that for any integer $n$ \[|\{\pi\in\D: \O(\pi)=n\}| = |\{\pi\in\P: |\pi|=n\}|.\] This result has been proven and re-discovered by many over the years. For example, P. Mork \cite{Mork} gave a solution to Schmidt's problem, the second author gave an independent proof in 2018 \cite{Uncu} without the knowledge of Schmidt's problem, and recently Andrews--Paule \cite{AndrewsPaule}. Andrews--Paule paper is particularly important because it led many new researchers into this area and many new proofs of Schmidt's theorem and its variants started appearing in a short period of time. Some examples that give novel proofs to Schmidt's problem are Alladi \cite{AlladiSchmidt}, Li-Yee \cite{LiYee}, and Bridges joint with the second author \cite{BridgesUncu}.

The original way the second author discovered and proved Schmidt's problem in 2016 uses Boulet--Stanley weights \cite{Boulet}. While working on her doctorate under the supervision of R. Stanley, Boulet wrote a 4-parameter generalization of the generating functions for the number of partitions where, given a Ferrers diagram, they decorated the odd indexed parts with alternating variables $a$ and $b$ (starting with $a$) and decorated the even indexed parts with alternating $c$ and $d$ (starting with $c$). In this setting, instead of $q$ keeping track of the size, we have a four variable function $\omega_\pi(a,b,c,d)$ where each variable counts the number of respective variables in the decorated Ferrers diagram. For example, the 4-decorated Ferrers diagram of the $\pi=(12,10,7,5,2)$ and the respective weight function are given next.

\begin{center}
\begin{tikzpicture}[line cap=round,line join=round,x=0.45cm,y=0.45cm]
\clip(0.5,0) rectangle (24.5,6.5);
\draw [line width=.25pt] (2,1)-- (4,1);
\draw [line width=.25pt] (2,2)-- (7,2);
\draw [line width=.25pt] (2,3)-- (9,3);
\draw [line width=.25pt] (2,4)-- (12,4);
\draw [line width=.25pt] (14,5)-- (2,5);
\draw [line width=.25pt] (2,6)-- (14,6);
\draw [line width=.25pt] (14,6)-- (14,5);
\draw [line width=.25pt] (13,6)-- (13,5);
\draw [line width=.25pt] (12,6)-- (12,4);
\draw [line width=.25pt] (11,6)-- (11,4);
\draw [line width=.25pt] (10,6)-- (10,4);
\draw [line width=.25pt] (9,6)-- (9,3);
\draw [line width=.25pt] (8,6)-- (8,3);
\draw [line width=.25pt] (7,2)-- (7,6);
\draw [line width=.25pt] (6,6)-- (6,2);
\draw [line width=.25pt] (5,6)-- (5,2);
\draw [line width=.25pt] (4,1)-- (4,6);
\draw [line width=.25pt] (3,6)-- (3,1);
\draw [line width=.25pt] (2,6)-- (2,1);
\draw (2.5,5.5) node[anchor=center] {$a$};
\draw (3.5,5.5) node[anchor=center] {$b$};
\draw (4.5,5.5) node[anchor=center] {$a$};
\draw (5.5,5.5) node[anchor=center] {$b$};
\draw (6.5,5.5) node[anchor=center] {$a$};
\draw (7.5,5.5) node[anchor=center] {$b$};
\draw (8.5,5.5) node[anchor=center] {$a$};
\draw (9.5,5.5) node[anchor=center] {$b$};
\draw (10.5,5.5) node[anchor=center] {$a$};
\draw (11.5,5.5) node[anchor=center] {$b$};
\draw (12.5,5.5) node[anchor=center] {$a$};
\draw (13.5,5.5) node[anchor=center] {$b$};

\draw (2.5,4.5) node[anchor=center] {$c$};
\draw (3.5,4.5) node[anchor=center] {$d$};
\draw (4.5,4.5) node[anchor=center] {$c$};
\draw (5.5,4.5) node[anchor=center] {$d$};
\draw (6.5,4.5) node[anchor=center] {$c$};
\draw (7.5,4.5) node[anchor=center] {$d$};
\draw (8.5,4.5) node[anchor=center] {$c$};
\draw (9.5,4.5) node[anchor=center] {$d$};
\draw (10.5,4.5) node[anchor=center] {$c$};
\draw (11.5,4.5) node[anchor=center] {$d$};

\draw (2.5,3.5) node[anchor=center] {$a$};
\draw (3.5,3.5) node[anchor=center] {$b$};
\draw (4.5,3.5) node[anchor=center] {$a$};
\draw (5.5,3.5) node[anchor=center] {$b$};
\draw (6.5,3.5) node[anchor=center] {$a$};
\draw (7.5,3.5) node[anchor=center] {$b$};
\draw (8.5,3.5) node[anchor=center] {$a$};

\draw (2.5,2.5) node[anchor=center] {$c$};
\draw (3.5,2.5) node[anchor=center] {$d$};
\draw (4.5,2.5) node[anchor=center] {$c$};
\draw (5.5,2.5) node[anchor=center] {$d$};
\draw (6.5,2.5) node[anchor=center] {$c$};

\draw (2.5,1.5) node[anchor=center] {$a$};
\draw (3.5,1.5) node[anchor=center] {$b$};

\draw (19,3.3) node[anchor=center] {,$\ \ \omega_\pi(a,b,c,d) = a^{11}b^{10}c^{8}d^7$.};
\end{tikzpicture}
\end{center} 

Boulet showed that the 4 parameter decorated generating functions for the ordinary partitions and for partitions into distinct parts also have product representations. 
\begin{theorem}[Boulet] \label{thm:boulet}For variables $a$, $b$, $c$, and $d$ and $Q:=abcd$, we have\begin{align} \label{eq:PSI}\Psi(a,b,c,d) &:= \sum_{\pi\in D} \omega_\pi(a,b,c,d) = \frac{(-a,-abc;Q)_\infty}{(ab;Q)_\infty},\\
\label{eq:PHI}\Phi(a,b,c,d) &:= \sum_{\pi\in U} \omega_\pi(a,b,c,d) = \frac{(-a,-abc;Q)_\infty}{(ab,ac,Q;Q)_\infty}.
\end{align}\end{theorem} 

It is easy to check that with the trivial choice $a=b=c=d=q$, $\omega_\pi(q,q,q,q)=q^{|\pi|}$, the 4-parameter generating functions become the generating function for the number of partitions and number of partitions into distinct parts, respectively. The choice $(a,b,c,d)=(q,q,1,1)$, proves Schmidt problem and yields another similar theorem of the sort for ordinary partitions:
\begin{theorem}[Theorems 1.3 and 6.3 \cite{Uncu}]
\begin{align}\label{Schmidt}\Psi(q,q,1,1) &= \sum_{\pi\in\D} q^{\O(\pi)} = \frac{1}{(q;q)_\infty},\\\label{UAP}\Phi(q,q,1,1) &= \sum_{\pi\in\P} q^{\O(\pi)} = \frac{1}{(q;q)^2_\infty}.\end{align}
\end{theorem}
This was the approach of the second author \cite{Uncu}. In the same paper, he also showed that Alladi's weighted theorem on Rogers-Ramanujan type partitions \cite{AlladiWeighted} is equivalent to Schmidt's problem.

Another independent proof of the right most equality of \eqref{UAP} recently appeared in the paper Andrews--Paule \cite[Theorem 2]{AndrewsPaule}, where they also proved Schmidt's problem. This also created an influx of research interests towards that one particular identity \eqref{UAP}. Some examples of these new works proving, dissecting and refining \eqref{UAP} are \cite{AlladiSchmidt, AndrewsKeith, BridgesUncu, ChernYee, Ji, LiYee}. In some of these texts, \eqref{UAP} is referred to as the Andrews--Paule theorem, however considering second author's earlier discovery and proof of this result \cite[Theorem 6.3]{Uncu}, we will refer to this result as the Uncu--Andrews--Paule theorem. 

In this paper, we present two refinements of the Uncu--Andrews--Paule theorem. The doubly refined theorem of this result is as follows. Let $\P_{\leq N}$ be the set of partitions with parts $\leq N$, and let \[\gamma(\pi) :=\lambda_1-\lambda_2+\lambda_3-\lambda_4+\dots\] be the \textit{alternating sum of parts} of a given partition $\pi=(\lambda_!,\lambda_2,\dots)$.

\begin{theorem}\label{thm:RefinedAlterPhi}Let $U_{N,j}(n)$ be the number of partitions $\pi\in\P_{\leq N}$ such that $\O(\pi)=n$ and $\gamma(\pi)= j$ and let $T_{N,j}(n)$ be the number of 2-color partitions (red and green)  of $n$ such that exactly $j$ parts are red and the largest green part does not exceed $N-j$. Then, \[U_{N,j}(n) = T_{N,j}(n).\]
\end{theorem}

We will also give various results involving the \textit{even-indexed parts' sum} statistics \[\E(\pi):= \lambda_2+\lambda_4+\lambda_6+\dots,\] where $\pi=(\lambda_1,\lambda_2,\lambda_3,\dots)$. One of such relations is as given below.

\begin{theorem}\label{cor:Intro} The number of partitions $\pi$ into distinct parts where $\E(\pi)=n$ and $\gamma(\pi)=j$ is equal to the number of partitions of $n$ into parts $\leq j$.
\end{theorem}

An example of this theorem with $n=5$ and $j=4$ is given in the next table.

\[\begin{array}{c|c}
\E(\pi)=5\ \& \ \gamma(\pi)=4 & |\pi|=5 \ \& \ \text{parts}\leq 4\\ \hline
(9,5),\ (8,5,1),\ (7,5,2)	& (4,1),\ (3,2),\  (3,1,1),\\
(7,4,2,1),\ (6,5,3),\ (6,4,3,1).	& (2,2,1),\ (2,1,1,1),\ (1,1,1,1,1).
\end{array}\]

\vspace{1mm}
The organization of this paper is as follows. In Section~\ref{Sec2}, we will recall some refinements of the generating functions proven by Boulet, compare it with a result by Berkovich and Warnaar on Rogers--Szeg\H{o} polynomials to get some weighted partition identities involving the alternating sum of parts and the sum of odd-indexed parts statistics. In Section~\ref{Sec3}, we will give finite analogues of Schmidt's problem and Uncu--Andrews--Paule theorem. Section~\ref{Sec4} is reserved for other Schmidt's problem like implications of the results in Section~\ref{Sec2}, which involve various alternating sums that are governed by the even indexed parts. Section~\ref{Sec5} has the Schmidt type results that focus on adding even indexed parts rather than the odd indexed parts. Finally, Section~\ref{Sec6} has the weighted counts connections between the Schmidt type statistics and multiplicative weights that were earlier studied by Alladi \cite{AlladiWeighted} and by the second author \cite{Uncu}.

\section{Refinements for Boulet's generating functions and Rogers--Szeg\H{o} polynomials}\label{Sec2}
%%%%%%%%%%%%%%%%%%%%%%%%%%%%%%%%%%%%%%%%%%%%%%%%%%%%%%%%%%%%%%%%%%%%%%%%%%%%%%%%%%%%%%%

In \cite{BU2} the authors made an extensive study on Boulet by imposing bounds on the largest part and the number of parts of the partitions. It should be noted that Ishikawa and Zeng \cite{Masao} were the first ones to present four variable generating functions for distinct and ordinary partitions by imposing a single bound on the largest part of partitions. In \cite{BU2}, the authors gave a different representation of the singly bounded generating functions and also gave two doubly bounded generating functions. Later, Fu and Zeng \cite{Fu} worked on the questions and techniques discussed in \cite{BU1, BU2}, and they presented doubly bounded generating functions with uniform bounds for the 4-decorated Ferrers diagrams of ordinary partitions and of partitions into distinct parts.

Let $\P_{\leq N}$ and $\D_{\leq N}$ be the sets of partitions from the sets $\P$ and $\D$, respectively, with the extra bound $N$ on the largest part. Define the generating functions \begin{align}
 \label{PSI_N}\Psi_N(a,b,c,d) &:= \sum_{\pi\in \D_{\leq N}} \omega_\pi(a,b,c,d),\\
 \label{PHI_N}\Phi_N(a,b,c,d) &:= \sum_{\pi\in \P_{\leq N}} \omega_\pi(a,b,c,d),
 \end{align}
which are finite analogues of Boulet's generating functions for the weighted count of 4-variable decorated Ferrers diagrams. In \cite{Masao}, Ishikawa and Zeng wrote explicit formulas for \eqref{PSI_N} and \eqref{PHI_N} using Pfaffians. 

\begin{theorem}[Ishikawa--Zeng]\label{MasaoThm} For a non-zero integer $N$, variables $a$, $b$, $c$, and $d$, we have
\begin{align}
\label{PSI2Nnu}\Psi_{2N+\nu}(a,b,c,d) &= \sum_{i=0}^N \genfrac{[}{]}{0pt}{}{N}{i}_{Q} (-a;Q)_{N-i+\nu}(-c;Q)_{i} (ab)^{i},\\
\label{PHI2Nnu}\Phi_{2N+\nu}(a,b,c,d) &=\frac{1}{(ac;Q)_{N+\nu}(Q;Q)_{N}}\sum_{i=0}^N \genfrac{[}{]}{0pt}{}{N}{i}_{Q} (-a;Q)_{N-i+\nu}(-c;Q)_{i} (ab)^{i},
\end{align}
where $\nu\in\{0,1\}$ and $Q=abcd$.
\end{theorem}

In Theorem~\ref{MasaoThm} and throughout the rest of the paper we use the standard definition \cite{theoryofpartitions} of the \textit{$q$-binomial coefficients}:
\[{n+m\brack n}_q := \left\{\begin{array}{cl}\displaystyle
\frac{(q;q)_{n+m}}{(q;q)_n(q;q)_m}, & \text{if}\ \ n,m\geq 0,\\[-1.5ex]\\
0, & \text{otherwise}.
\end{array}\right.\]

In \cite{BU1}, we gave a companion to Theorem~\ref{MasaoThm}.

\begin{theorem}\label{finiteBoulet} For a non-zero integer $N$, variables $a$, $b$, $c$, $d$, and $Q=abcd$, we have
\begin{align}
\label{finiteBouletPSI}\Psi_{2N+\nu}(a,b,c,d) &=\sum_{i=0}^N \genfrac{[}{]}{0pt}{}{N}{i}_{Q} (-a;Q)_{i+\nu}(-abc;Q)_{i} \frac{(ac;Q)_{N+\nu}}{(ac;Q)_{i+\nu}}(ab)^{N-i},\\
\label{finiteBouletPHI}\Phi_{2N+\nu}(a,b,c,d) &=\frac{1}{(Q;Q)_{N}}\sum_{i=0}^N \genfrac{[}{]}{0pt}{}{N}{i}_{Q} \frac{(-a;Q)_{i+\nu}(-abc;Q)_{i} }{(ac;Q)_{i+\nu}}(ab)^{N-i},
\end{align} where $\nu\in\{0,1\}$.
\end{theorem}
It should be noted that our derivation of Theorem~\ref{finiteBoulet}, unlike Theorem~\ref{MasaoThm}, is completely combinatorial. 

Another useful connection is the first author and Warnaar's theorem \cite{WarnaarBerkovich} on the Rogers--Szeg\H{o} polynomials \[H_N(z,q):=\sum_{k=0}^N \genfrac{[}{]}{0pt}{}{N}{k}_{q}z^k.\]

\begin{theorem} [Berkovich--Warnaar]\label{HermiteBerkovichTHM} Let $N$ be a non-negative integer, then the Rogers--Szeg\H{o} polynomials can be expressed as
\begin{equation}\label{HermiteBerkovich}
H_{2N+\nu}(zq,q^2) = \sum_{k=0}^{N}\genfrac{[}{]}{0pt}{}{N}{k}_{q^4} (-zq;q^4)_{N-k+\nu}(-q/z;q^4)_k (zq)^{2k},
\end{equation} where $\nu=0$ or $1$.
\end{theorem}

We see that the right-hand sides of \eqref{PSI2Nnu} and \eqref{HermiteBerkovich} coincide for a particular choice of $(a,b,c,d)$:
\begin{equation}\label{RogersSzegoToPSI}H_{N}(zq,q^2) = \Psi_{N}(zq,zq,q/z,q/z).\end{equation}
In \cite{BU2}, authors showed that \begin{equation}\label{eq:511}\Psi_{N}(zq,zq,q/z,q/z) = \sum_{\pi\in\D_{\leq N}} q^{|\pi|} z^{\gamma(\pi)},\end{equation} where $\gamma(\pi):=\lambda_1-\lambda_2+\lambda_3-\lambda_4+\dots,$ the \textit{alternating sum of the parts} of $\pi=(\lambda_1,\lambda_2,\lambda_3,\lambda_4,\dots)$. Moreover, in the same paper the authors, using \eqref{RogersSzegoToPSI}, also showed that the coefficient of the term $z^k$ in \eqref{eq:511} is $q^k{N\brack k}_{q^2}$. 

Next, by first shifting $z$ to $zq$ and then mapping $q^2\mapsto q$, yields

\begin{theorem}\label{thm:SchmidtTypeAlternating}
\begin{align}
\label{eq:Psi_alternating}\Psi_N(qz,qz,1/z,1/z) &= \sum_{k = 0}^N q^kz^k {N\brack k}_q ,\\
\label{eq:Phi_alternating}\Phi_N(qz,qz,1/z,1/z) &= \frac{\Psi_N(qz,qz,1/z,1/z)}{(q;q)_N} = \sum_{k = 0}^N \frac{q^k z^k}{(q;q)_k(q;q)_{N-k}}.
\end{align}
\end{theorem}

One can rewrite the above theorem as

\begin{theorem}\label{thm:SchmidtTypeAlternatingComb}
\begin{align}
\label{eq:Psi_alternatingComb}\sum_{\pi\in\D_{\leq N}} q^{\O(\pi)} z^{\gamma(\pi)} &= \sum_{k = 0}^N q^k z^k{N\brack k}_q,\\
\label{eq:Phi_alternatingComb}\sum_{\pi\in\P_{\leq N}} q^{\O(\pi)} z^{\gamma(\pi)}&= \sum_{k = 0}^N \frac{q^k z^k}{(q;q)_k(q;q)_{N-k}}.
\end{align}
\end{theorem}

\section{Another look at the refined Schmidt and results alike}\label{Sec3}

Theorem~\ref{MasaoThm} (and Theorem~\ref{finiteBoulet}) gives a sum representation for the bounded analogues of Schmidt's problem \eqref{Schmidt} and the Uncu--Andrews--Paule theorem \eqref{UAP}. 

\begin{theorem} Let $N$ be any integer and $\nu =0$ or 1, then
\begin{align}
\label{eq:Old_PsiN_rep}\Psi_{2N+\nu}(q,q,1,1) &=\sum_{\pi\in\D_{\leq 2N+\nu}} q^{\O(\pi)}= \sum_{i=0}^N \genfrac{[}{]}{0pt}{}{N}{i}_{q^2} (-q;q^2)_{N-i+\nu}(-1;q^2)_{i} q^{2i},\\
\label{eq:Old_PhiN_rep}\Phi_{2N+\nu}(q,q,1,1) &=\sum_{\pi\in\P_{\leq 2N+\nu}}q^{\O(\pi)}=\frac{1}{(q;q)_{2N+\nu}}\sum_{i=0}^N \genfrac{[}{]}{0pt}{}{N}{i}_{q^2} (-q;q^2)_{N-i+\nu}(-1;q^2)_{i} q^{2i}.
\end{align}
\end{theorem}

One can also prove that these generating functions have an alternate representations.

\begin{theorem}\label{thm:colors} Let $N$ be any integer then
\begin{align}
\label{eq:New_PsiN_rep}\sum_{\pi\in\D_{\leq N}} q^{\O(\pi)}&= \sum_{i=0}^{N} \genfrac{[}{]}{0pt}{}{N}{i}_{q} q^{i},\\
\label{eq:New_PhiN_rep}\sum_{\pi\in\P_{\leq N}}q^{\O(\pi)}&=\sum_{i=0}^{N} \frac{q^i}{(q;q)_i(q;q)_{N-i}}.
\end{align}
\end{theorem}

Picking $z=1$ in Theorem~\ref{thm:SchmidtTypeAlternatingComb} gives Theorem~\ref{thm:colors}.

We can also give new partition interpretations for the sums in \eqref{eq:New_PsiN_rep} and \eqref{eq:New_PhiN_rep}. This directly leads to new weighted partition identities once we compare them with the previously known interpretations of $\Psi_N(q,q,1,1)$ and $\Phi_N(q,q,1,1)$. We start with the new weighted partition identity related to the Schmidt's problem. 

\begin{theorem}\label{thm:FinWeightedSchmidt} Let $S_N(n)$ be the number of partitions into distinct parts $\pi = (\lambda_1,\lambda_2,\dots)$ with $\O(\pi)=n$ and $\lambda_1\leq N$. Let $\Gamma_N(n)$ be the number of partitions of $n$ where the largest hook length is $\leq N$. Then, \[S_N(n) = \Gamma_N(n).\]
\end{theorem}

Recall that, on a Ferrers diagram, a \textit{hook length} of a box is defined as one plus the number of boxes directly to the right and directly below to the chosen box. It is then easy to see that the largest hook length is the hook length of the top-left-most box in a Ferrers diagram. Whence, we could also define $\Gamma_N(n)$ as the number of partitions of $n$ where the quantity ``the number of parts plus the largest part minus 1" is less than or equal to $N$.

We exemplify Theorem~\ref{thm:FinWeightedSchmidt} with $n=N=4$ and list the partitions counted by $S_N(n)$ and $\Gamma_N(n)$ here.

\[\begin{array}{c|c}
S_4(4)=5 &  \Gamma_4(4)=5\\ \hline
(4,3)	&	(4)			\\
(4,2)	&	(3,1)		\\
(4,1)	&	(2,2)		\\
(4)		&	(2,1,1)		\\
(3,2,1)	&	(1,1,1,1)	\\	
\end{array}
\]

Notice that, whereas most partitions in this example have their largest hook lengths to be the upper limit 4, the partition $(2,2)$ has the largest hook length $3<4$.

The proof of Theorem~\ref{thm:FinWeightedSchmidt} only requires us to interpret \eqref{eq:New_PsiN_rep} as the generating function for the $\Gamma_N(n)$ numbers. To that end, we focus on the summand of the right-hand side in \eqref{eq:New_PsiN_rep}. Recall that the $q$-binomial coefficient ${N\brack i}_q$ is the generating function of partitions that fit in a box with $i$ rows and $(N-i)$-columns. We fit a column of height $i$ in front of that box. There are exactly $i$ boxes on the left most column of this Ferrers diagram and up to $N-i$ boxes on the top most row. Hence, the largest hook length of this partition is $\leq N$. By summing over all possible $i$, we cover all possible column heights. This shows that the right-hand side of \eqref{eq:New_PsiN_rep} is the generating function of $\Gamma_N(n)$ where $q$ is keeping track of the size of the partitions.

A similar weighted theorem comes from comparing the known interpretation of $\Phi_N(q,q,1,1)$ and the new interpretation of the same object coming from \eqref{eq:New_PhiN_rep}.

\begin{theorem} Let $U_N(n)$ be the number of partitions  $\pi = (\lambda_1,\lambda_2,\dots)$ with $\O(\pi)=n$ and $\lambda_1\leq N$. Let $T_N(n)$ be the number of 2-color partitions (red and green)  of $n$ such that the number of parts in red plus the size of the largest green part does not exceed $N$. Then, \[U_N(n) = T_N(n).\]
\end{theorem}

We give an example of this theorem with $n=4$ and $N=3$. We will use subscripts $r$ and $g$ to indicate colors of the parts of partitions while listing the partitions related to $T_3(4).$

\[
\begin{array}{c|c}
U_3(4) = 15 & T_3(4) = 15\\
\hline
\begin{array}{cc} 
(3,3,1,1), & 	(3,3,1),\\
(3,2,1,1), &	(3,2,1),\\
(3,1,1,1), &	(3,1,1),\\
(2,2,2,2),	&	(2,2,2),\\
(2,2,1,1,1,1),&	(2,2,1,1,1),\\
(2,2,2,1)	&	(2,1,1,1,1,1),\\
(2,1,1,1,1),&	(1,1,1,1,1,1,1),\\
\end{array} & 
\begin{array}{cc} 
(4_{\color{red} r} )	&	(3_{\color{red} r} ,1_{\color{red} r} ) \\
(3_{\color{red} r} ,1_{\color{olive} g} ) & (3_{\color{olive} g} ,1_{\color{olive} g} ) \\
(2_{\color{red} r} ,2_{\color{olive} g} ),& (2_{\color{olive} g} ,2_{\color{olive} g} ), \\
(2_{\color{red} r} ,1_{\color{red} r} ,1_{\color{red} r} ), &(2_{\color{red} r} ,1_{\color{red} r} ,1_{\color{olive} g} ), \\
(2_{\color{red} r} ,1_{\color{olive} g} ,1_{\color{olive} g} ), &(2_{\color{olive} g} ,1_{\color{red} r} ,1_{\color{olive} g} ), \\
(2_{\color{olive} g} ,1_{\color{olive} g} ,1_{\color{olive} g} ), &(1_{\color{red} r} ,1_{\color{red} r} ,1_{\color{red} r} ,1_{\color{olive} g} )\\
(1_{\color{red} r} ,1_{\color{red} r} ,1_{\color{olive} g} ,1_{\color{olive} g} ),& (1_{\color{red} r} ,1_{\color{olive} g} ,1_{\color{olive} g} ,1_{\color{olive} g} ),\\
\end{array}\\
(1,1,1,1,1,1,1,1). & (1_{\color{olive} g} ,1_{\color{olive} g} ,1_{\color{olive} g} ,1_{\color{olive} g} ).
\end{array}
\]

Note that the count of $T_N(n)$ is not symmetric in colors red and green. In our example, partitions such as $(4_g)$, $(2_g,1_r,1_r)$ and $(1_r,1_r,1_r,1_r)$ are not counted.

The equation \eqref{eq:Psi_alternatingComb} can also be presented as a refined Schmidt problem like result.

\begin{theorem}\label{thm:RefinedAlterPsi}Let $S_{N,j}(n)$ be the number of partitions $\pi\in\D_{\leq N}$ such that $\O(\pi)=n$ and $\gamma(\pi) = j$ and let $\Gamma_{N,j}(n)$ be the number of partitions into exactly $j$ parts where the largest hook length is $\leq N$. Then, \[S_{N,j}(n) = \Gamma_{N,j}(n).\]
\end{theorem}

The proof of this result comes from comparing and interpreting the coefficients of $z^j$ in \eqref{eq:Psi_alternatingComb}. Similarly, we can get an analogous result for ordinary partitions by comparing the $z^j$ terms in \eqref{eq:Phi_alternatingComb}, which yields Theorem~\ref{thm:RefinedAlterPhi}. Moreover, letting $N$ tend to infinity we arrive at \cite[Corollary 2]{AndrewsKeith}. This corollary was implicit in Mork's original solution \cite{Mork} of the Schmidt problem. 

We give examples of Theorems \ref{thm:RefinedAlterPsi} and \ref{thm:RefinedAlterPhi} in the next table.

\[
\begin{array}{cc|cc}
S_{3,2}(4) = 2  &   \Gamma_{3,2}(4) = 2 &  U_{3,1}(4) = 6  &   T_{3,1}(4) =  6\\ \hline
	(4,2)	&	(3,1)	&	(3,3,1)			&	(4_{\color{red} r})	\\
	(3,2,1)	&	(2,2)	&	(3,2,1,1)		&	(3_{\color{red} r},1_{\color{olive} g})	\\
		&		&	(2,2,1,1,1)		&	(2_{\color{red} r},1_{\color{olive} g},1_{\color{olive} g})	\\
		&		&	(2,2,2,1)		&	(2_{\color{red} r},2_{\color{olive} g})	\\
		&		&	(2,1,1,1,1,1)	&	(1_{\color{red} r},1_{\color{olive} g},1_{\color{olive} g},1_{\color{olive} g})	\\
		&		&	(1,1,1,1,1,1,1)	&	(2_{\color{olive} g},1_{\color{red} r},1_{\color{olive} g})	\\
\end{array}
\]

\section{Schmidt type results}\label{Sec4}

For a given partition $\pi$, choices of variables in Boulet's and subsequently the bounded analogues of Boulet's theorem can give us great insight into Schmidt type partition theorems. For example, by the choice $(a,b,c,d) = (q,q,-1,-1)$ in \eqref{eq:PSI}, we directly get the following theorem.

\begin{theorem}\label{thm:distinct1mod2}
\begin{align}
\label{eq:distinct1mod2} \sum_{\pi\in\D} (-1)^{\E(\pi)} q^{\O(\pi)} &= (-q;q^2)_\infty,\\
\label{eq:ordinary0mod2} \sum_{\pi\in\P} (-1)^{\E(\pi)} q^{\O(\pi)} &= \frac{1}{(q^2;q^2)_\infty}.
\end{align}
\end{theorem}

Theorem~\ref{thm:distinct1mod2} translates into the following two weighted combinatorial result.

\begin{theorem}
The number of partitions $\pi$ into distinct parts counted with weight $(-1)^{\E(\pi)}$ such that $\O(\pi)=n$ is equal to the number of partitions of $n$ into distinct odd parts.
\end{theorem}

\begin{theorem}
The number of partitions $\pi$ counted with weight $(-1)^{\E(\pi)}$ such that $\O(\pi)=n$  is equal to the number of partitions of $n$ into even parts.\end{theorem}

The finite analogues of Boulet's theorem can directly be used to get refinements of these results.

\begin{theorem}\label{thm:refined1} Let $\D_{\leq N}$ be the set of partitions into distinct parts $\leq N$, then
\begin{align}
\label{eq:distinct1mod2_N}  \sum_{\pi\in\D_{\leq N}} (-1)^{\E(\pi)} q^{\O(\pi)} = (-q;q^2)_{\lceil N/2\rceil},\\
\label{eq:ordinary0mod2_N}  \sum_{\pi\in\P_{\leq N}} (-1)^{\E(\pi)} q^{\O(\pi)} = \frac{1}{(q^2;q^2)_{\lceil N/2\rceil}},
\end{align} where $\lceil\cdot\rceil$ is the standard ceil function.
\end{theorem}

With the choice $c=-1$, the sums in Theorem~\ref{MasaoThm} collapse to the $i=0$ term only and this yields a direct proof of Theorem~\ref{thm:refined1}.

For example let $N=2$, then the explicit weighted count of the partitions related to \eqref{eq:distinct1mod2_N} is as follows:
\[
\begin{array}{c}
\begin{array}{cc|cc|cc}
\pi\in\D_{\leq4} & (-1)^{\E(\pi)} q^{\O(\pi)} & \pi\in\D_{\leq4} & (-1)^{\E(\pi)} q^{\O(\pi)} & \pi\in\D_{\leq4} & (-1)^{\E(\pi)} q^{\O(\pi)}\\
\hline
 (4)	& q^4	& (4,2,1)	& q^5	& (2)	& q^2	\\
 (4,3)	& -q^4	& (4,3,2,1)	& q^6	& (2,1)	& -q^2 	\\
 (4,2)	& q^4	& (3)		& q^3	& (1)	& q	\\
 (4,1)	& -q^4	& (3,2)		& q^3	& \emptyset	& 1	\\
(4,3,2) & -q^6	& (3,1)		& -q^3	& 	& 	\\
(4,3,1)	& -q^5	& (3,2,1)	& q^4	& 	& 	
\end{array}\\[-1.5ex]\\
\displaystyle
\Psi_{4}(q,q,-1,-1) =\sum_{\pi\in\D_{\leq 4}} (-1)^{\E(\pi)} q^{\O(\pi)}= q^4+q^3+q+1 = (1+q)(1+q^3) = (-q;q^2)_2.
\end{array}
\]

We can also choose $(a,b,c,d)= (q,q,-1,1)$ and get another theorem.

\begin{theorem}\label{thm:distinct1mod2Neg}
\begin{align}
\label{eq:distinctNeg1mod2} \sum_{\pi\in\D} (-1)^{\lceil\E\rceil(\pi)} q^{\O(\pi)} &= (-q;-q^2)_\infty = (-q,q^3;q^4)_\infty,\\
\label{eq:ordinaryNeg0mod2} \sum_{\pi\in\P} (-1)^{\lceil\E\rceil(\pi)} q^{\O(\pi)} &= \frac{1}{(-q^2;-q^2)_\infty} = \frac{1}{(-q^2,q^4;q^4)_\infty},
\end{align} where $\lceil\E\rceil (\lambda_1,\lambda_2,\lambda_3,\dots) = \lceil\lambda_2/2\rceil+\lceil\lambda_4/2\rceil+\lceil\lambda_6/2\rceil+\dots$.
\end{theorem}

The refinement of Theorem~\ref{thm:distinct1mod2Neg} can also be found in a similar manner.

\begin{theorem}\label{thm:refined2}
\begin{align}
\label{eq:distinctNeg1mod2_N}  \sum_{\pi\in\D_{\leq N}} (-1)^{\lceil\E(\pi)\rceil} q^{\O(\pi)} &= (-q;-q^2)_{\lceil N/2\rceil},\\
\label{eq:ordinaryNeg0mod2_N} \sum_{\pi\in\P_{\leq N}} (-1)^{\lceil\E\rceil(\pi)} q^{\O(\pi)} &= \frac{1}{(-q^2;-q^2)_{\lceil N/2\rceil}}.
\end{align}
\end{theorem}

For $N=2$ \eqref{eq:distinctNeg1mod2_N} is as follows:
\[
\begin{array}{c}
\begin{array}{cc|cc|cc}
\pi\in\D_{\leq4} & (-1)^{\lceil\E(\pi)\rceil} q^{\O(\pi)} & \pi\in\D_{\leq4} & (-1)^{\lceil\E(\pi)\rceil} q^{\O(\pi)} & \pi\in\D_{\leq4} & (-1)^{\lceil\E(\pi)\rceil} q^{\O(\pi)}\\
\hline
 (4)	& q^4	& (4,2,1)	& -q^5	& (2)	& q^2	\\
 (4,3)	& q^4	& (4,3,2,1)	& -q^6	& (2,1)	& -q^2 	\\
 (4,2)	& -q^4	& (3)		& q^3	& (1)	& q	\\
 (4,1)	& -q^4	& (3,2)		& -q^3	& \emptyset	& 1	\\
(4,3,2) & q^6	& (3,1)		& -q^3	& 	& 	\\
(4,3,1)	& q^5	& (3,2,1)	& -q^4	& 	& 	
\end{array}\\[-1.5ex]\\
\displaystyle
\Psi_{4}(q,q,-1,1) =\sum_{\pi\in\D_{\leq 4}} (-1)^{\lceil\E(\pi)\rceil} q^{\O(\pi)}= -q^4-q^3+q+1 = (1+q)(1-q^3) = (-q;q^2)_2.
\end{array}
\]

A modulo 3 related group of Schmidt type results come from the following choice $(a,b,c,d)= (q,q,-1,-q)$.

\begin{theorem}\label{thm:distinct1mod2NegFloor}
\begin{align}
\label{eq:distinctNeg1mod2Floor} \sum_{\pi\in\D} (-1)^{\E(\pi)} q^{\O(\pi)+\lfloor \E \rfloor (\pi)} &=(-q;q^3)_\infty,\\
\label{eq:ordinaryNeg0mod2Floor} \sum_{\pi\in\P} (-1)^{\E(\pi)} q^{\O(\pi)+\lfloor \E \rfloor (\pi)} &=\frac{1}{(q^3;q^3)_\infty},
\end{align} where $\lfloor\E\rfloor (\lambda_1,\lambda_2,\lambda_3,\dots) = \lfloor\lambda_2/2\rfloor+\lfloor\lambda_4/2\rfloor+\lfloor\lambda_6/2\rfloor+\dots$.
\end{theorem}

Once again, a refinement of Theorem~\ref{thm:distinct1mod2NegFloor} similar to the refinements above can easily be found.

\begin{theorem}\label{thm:refined3}
\begin{align}
\label{eq:distinctNeg1mod2_fN}  \sum_{\pi\in\D_{\leq N}} (-1)^{\lfloor\E(\pi)\rfloor} q^{\O(\pi)} &= (-q;q^3)_{\lceil N/2\rceil},\\
\label{eq:ordinaryNeg0mod2_fN} \sum_{\pi\in\P_{\leq N}} (-1)^{\lfloor\E\rfloor(\pi)} q^{\O(\pi)} &= \frac{1}{(q^3;q^3)_{\lceil N/2\rceil}}.
\end{align}
\end{theorem}

We add one explicit example for \eqref{eq:distinctNeg1mod2_fN} with $N=2$:
\[
\begin{array}{c}
\begin{array}{cc|cc|cc}
\pi\in\D_{\leq4} & (-1)^{\lfloor\E(\pi)\rfloor} q^{\O(\pi)} & \pi\in\D_{\leq4} & (-1)^{\lfloor\E(\pi)\rfloor} q^{\O(\pi)} & \pi\in\D_{\leq4} & (-1)^{\lfloor\E(\pi)\rfloor} q^{\O(\pi)}\\
\hline
 (4)	& q^4	& (4,2,1)	& q^6	& (2)	& q^2	\\
 (4,3)	& -q^5	& (4,3,2,1)	& q^7	& (2,1)	& -q^2 	\\
 (4,2)	& q^5	& (3)		& q^3	& (1)	& q	\\
 (4,1)	& -q^4	& (3,2)		& q^4	& \emptyset	& 1	\\
(4,3,2) & -q^7	& (3,1)		& -q^3	& 	& 	\\
(4,3,1)	& -q^6	& (3,2,1)	& q^5	& 	& 	
\end{array}\\[-1.5ex]\\
\displaystyle
\Psi_{4}(q,q,-1,-q) =\sum_{\pi\in\D_{\leq 4}} (-1)^{\lfloor\E(\pi)\rfloor} q^{\O(\pi)}= q^5+q^4+q+1 = (1+q)(1+q^4) = (-q;q^3)_2.
\end{array}
\]

In general, $c=-1$ choice in Theorem~\ref{thm:boulet} always leads to cancellation of  terms in the products \eqref{eq:PSI} and \eqref{eq:PHI}.

\begin{theorem}\label{thm:MainSec4}
Let $(a,b,c,d) = (q^r,q^t,-1,\varepsilon q^s)$ for some non-negative integers $r,\, t,\, s$, (with $t+r>0$) and $\varepsilon=\pm 1$. Then we have
\begin{align}
\label{eq:distinctGeneral} \sum_{\pi\in\D} (-1)^{\lceil\E\rceil(\pi)} \varepsilon^{\lfloor\E\rfloor(\pi)} q^{r\lceil\O\rceil(\pi)+t\lfloor \O \rfloor (\pi)+ s \lfloor \E \rfloor (\pi)} &=(-q^r;\varepsilon  q^{r+t+s})_\infty,\\
\label{eq:ordinaryGeneral} \sum_{\pi\in\P} (-1)^{\lceil\E\rceil(\pi)} \varepsilon^{\lfloor\E\rfloor(\pi)} q^{r\lceil\O\rceil(\pi)+t\lfloor \O \rfloor (\pi)+ s \lfloor \E \rfloor (\pi)} &=\frac{1}{(\varepsilon  q^{r+t+s};\varepsilon q^{r+t+s})_\infty},
\end{align}
where $\lceil\O\rceil(\pi)$ and $\lfloor \O \rfloor(\pi) $ are defined in a similar fashion to $\lceil\E\rceil(\pi)$ and $\lfloor \E \rfloor(\pi)$.
\end{theorem}

Similarly the refinement of Theorem~\ref{thm:MainSec4} is as follows.

\begin{theorem}\label{thm:MainSec4N}
Let $(a,b,c,d) = (q^r,q^t,-1,\varepsilon q^s)$ for some non-negative integers $r,\, t,\, s$, (with $t+r>0$) and $\varepsilon=\pm 1$. Then we have
\begin{align}
\label{eq:distinctGeneral} \sum_{\pi\in\D_{\leq N}} (-1)^{\lceil\E\rceil(\pi)} \varepsilon^{\lfloor\E\rfloor(\pi)} q^{r\lceil\O\rceil(\pi)+t\lfloor \O \rfloor (\pi)+ s \lfloor \E \rfloor (\pi)} &=(-q^r;\varepsilon  q^{r+t+s})_{\lceil N/2\rceil},\\
\label{eq:ordinaryGeneral} \sum_{\pi\in\P_{\leq N}} (-1)^{\lceil\E\rceil(\pi)} \varepsilon^{\lfloor\E\rfloor(\pi)} q^{r\lceil\O\rceil(\pi)+t\lfloor \O \rfloor (\pi)+ s \lfloor \E \rfloor (\pi)} &=\frac{1}{(\varepsilon  q^{r+t+s};\varepsilon q^{r+t+s})_{\lceil N/2\rceil}}.
\end{align}
\end{theorem}

Theorems~\ref{thm:distinct1mod2}, \ref{thm:distinct1mod2Neg} and \ref{thm:distinct1mod2NegFloor} are special cases of Theorem~\ref{thm:MainSec4}. Similarly, Theorems~\ref{thm:refined1}, \ref{thm:refined2} and \ref{thm:refined3} are special cases of Theorem~\ref{thm:MainSec4N}. 

\section{Partitions identities for uncounted odd-indexed parts}\label{Sec5}

In sections~\ref{Sec3} and \ref{Sec4}, we were making sure that our substitutions to Boulet weights were never $a=b=1$ or $a=b=-1$. These choices clearly introduce a pole to the product representations of the generating functions \eqref{eq:PSI} and \eqref{eq:PHI}. However, these substitutions can be entertained in the finite analogue of these generating functions.

For example, substituting $a=b=-1$ and $c=d=q$ in \eqref{PSI2Nnu} and \eqref{PHI2Nnu} yields the two equations presented in the following theorem, respectively.

\begin{theorem} For some non-negative integer $N$,
\begin{align}
\label{eq:OddsMinus1Psi}\sum_{\pi\in\D_{\leq N}}(-1)^{\O(\pi)} q^{\E(\pi)} &= \left\{ \begin{array}{cl}
(-q;q^2)_{N/2} &\text{if }N\text{ is even},\\[-1.5ex]\\
0 &\text{if }N\text{ is odd}.
\end{array}\right.\\
\label{eq:OddsMinus1Phi}\sum_{\pi\in\P_{\leq N}} (-1)^{\O(\pi)} q^{\E(\pi)} &= \left\{ \begin{array}{cl}\displaystyle
\frac{1}{(q^2;q^2)_{N/2}} &\text{if }N\text{ is even},\\
0 &\text{if }N\text{ is odd}.
\end{array}\right.
\end{align}
\end{theorem}

%As one can see that neither \eqref{eq:OddsMinus1Psi}, nor \eqref{eq:OddsMinus1Phi} have a (uniform) limit as $N$ tends to infinity. (Although in each case there are two subsequences with clear limits.)

We list partitions and their respective weights below to demonstrate \eqref{eq:OddsMinus1Psi} with $N=4$ and $N=3$.

\[
\begin{array}{c}
\begin{array}{cc|cc|cc}
\pi\in\D_{\leq4} & (-1)^{\O(\pi)} q^{\E(\pi)} & \pi\in\D_{\leq4} & (-1)^{\O(\pi)} q^{\E(\pi)} & \pi\in\D_{\leq4} & (-1)^{\O(\pi)} q^{\E(\pi)}\\
\hline
 (4)	& 1		& (4,2,1)	& -q^2	& (2)	& 1	\\
 (4,3)	& q^3	& (4,3,2,1)	& q^4	& (2,1)	& q 	\\
 (4,2)	& q^2	& (3)		& -1	& (1)	& -1	\\
 (4,1)	& q		& (3,2)		& -q^2	& \emptyset	& 1	\\
(4,3,2) & q^3	& (3,1)		& -q	& 	& 	\\
(4,3,1)	& -q^3	& (3,2,1)	& q^2	& 	& 	
\end{array}\\[-1.5ex]\\
\displaystyle
\Psi_{4}(-1,-1,q,q) =\sum_{\pi\in\D_{\leq 4}} (-1)^{\O(\pi)} q^{\E(\pi)}= q^4+q^3+q+1 = (1+q)(1+q^3) = (-q;q^2)_2.\\[-1.5ex]\\
\Psi_{3}(-1,-1,q,q) =\sum_{\pi\in\D_{\leq 3}} (-1)^{\O(\pi)} q^{\E(\pi)}=0.
\end{array}
\]

We let $z=1/q$ and then map $q^2\mapsto q$ in \eqref{RogersSzegoToPSI} to get the next theorem.

\begin{theorem}\label{thm:Odds1}
\begin{align}
\label{eq:Odds1Psi}\sum_{\pi\in\D_{\leq N}} q^{\E(\pi)} &= \sum_{i=0}^N {N\brack i }_q,\\
\label{eq:Odds1Phi}\sum_{\pi\in\P_{\leq N}} q^{\E(\pi)} &= \sum_{i=0}^N \frac{1}{(q;q)_i(q;q)_{N-i}}.
\end{align}
\end{theorem}

We demonstrate this theorem with $N=4$ below.

\[
\begin{array}{cc|cc|cc}
\pi\in\D_{\leq4} &  q^{\E(\pi)} & \pi\in\D_{\leq4} &  q^{\E(\pi)} & \pi\in\D_{\leq4} &  q^{\E(\pi)}\\
\hline
 (4)	& 1		& (4,2,1)	& q^2	& (2)	& 1	\\
 (4,3)	& q^3	& (4,3,2,1)	& q^4	& (2,1)	& q 	\\
 (4,2)	& q^2	& (3)		& 1	& (1)	& 1	\\
 (4,1)	& q		& (3,2)		& q^2	& \emptyset	& 1	\\
(4,3,2) & q^3	& (3,1)		& q	& 	& 	\\
(4,3,1)	& q^3	& (3,2,1)	& q^2	& 	& 	
\end{array}
\]
\begin{align*}
\Psi_{4}(1,1,q,q) &=\sum_{\pi\in\D_{\leq 4}} q^{\E(\pi)}= 5+3q+4q^2+3q^3+q^4\\
&= 1+ (1+q+q^2+q^3) + (1+q+2q^2+q^3+q^4) + (1+q+q^2+q^3) +1 \\
&= {4\brack 0}_q + {4\brack 1}_q + {4\brack 2}_q+ {4\brack 3}_q +{4\brack 4}_q.
\end{align*}

This result can be refined if we keep track of the variable $z$ in \eqref{RogersSzegoToPSI}. Let $z\mapsto z/q$ and $q^2\mapsto q$ in this order in \eqref{RogersSzegoToPSI} to see two variable generalizations of \eqref{eq:Odds1Psi} and \eqref{eq:Odds1Phi}. Then, by extracting the coefficient of $z^j$, we get the following theorem.

\begin{theorem}\label{thm:RefinedOdds1} For a partition $\pi = (\lambda_1,\lambda_2,\dots)$, let $\gamma(\pi) = \lambda_1-\lambda_2+\lambda_3-\lambda_4+\dots$ be the \textit{alternating sum of the parts of }$\pi$. Then, 
\begin{align}
\label{eq:RefOdds1Psi}\sum_{\substack{\pi\in\D_{\leq N}\\\gamma(\pi) = j}} q^{\E(\pi)} &=  {N\brack j }_q,\\
\label{eq:RefOdds1Phi}\sum_{\substack{\pi\in\P_{\leq N}\\\gamma(\pi) = j}}q^{\E(\pi)} &= \frac{1}{(q;q)_j(q;q)_{N-j}}.
\end{align}
\end{theorem}

We also note a direct combinatorial proof for \eqref{eq:RefOdds1Psi} using the Sylvester's bijection \cite{Bressoud}. We know that ${N\brack j}_q$ is the generating function for the number of partitions with parts $\leq N-j$ and the number of parts $\leq j$. Let $\pi$ be one of such partitions into $\leq j$ parts where its largest part is $\leq N-j$. Next, we construct a partition $\pi_o$ into $j$ odd parts by first putting $j$ boxes in a column followed by putting the Ferrers diagram of the partition $\pi$ on the right-side of the $j$-box column, and its reflection \reflectbox{$\pi$} on to the left of the column. When it is read line-by-line, this is a partition into odd parts, and $|\pi_o| = 2|\pi| +j$. This construction from $\pi$ to $\pi_o$ is bijective and can be reversed easily. Now we apply the Sylvester's bijection to $\pi_o$. This takes $\pi_o$ to a partition into distinct parts $\pi_d$. Obviously the composition of the bijections $\pi\mapsto\pi_o\mapsto\pi_d$ is a bijection. The largest part of $\pi_d$ is the sum of the main column length (exactly $j$) and the size of the largest part of $\pi$ (which is $\leq N-j$), hence $\leq N$. It is clear that $\gamma(\pi_d) = j$. Finally, $\E(\pi_d)$ is equal to the number of boxes in the Ferrers diagram of the reflected copy in the Ferrers diagram, which is exactly $|\pi|$. One example of this construction is given in Figure~\ref{figure_ferrers}.

\definecolor{qqwuqq}{rgb}{0,0.39215686274509803,0}
\definecolor{ccqqqq}{rgb}{0.8,0,0}
\begin{center}
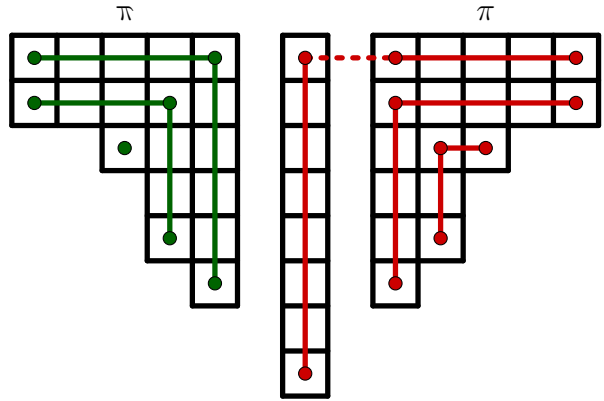
\begin{figure}[h]
\begin{tikzpicture}[line cap=round,line join=round,x=0.6cm,y=0.6cm]
\clip(0,0.8) rectangle (15.5,10);
\draw [line width=2pt] (1,9)-- (6,9);
\draw [line width=2pt] (1,8)-- (6,8);
\draw [line width=2pt] (1,7)-- (6,7);
\draw [line width=2pt] (3,6)-- (6,6);
\draw [line width=2pt] (4,5)-- (6,5);
\draw [line width=2pt] (4,4)-- (6,4);
\draw [line width=2pt] (5,3)-- (6,3);
\draw [line width=2pt] (9,3)-- (10,3);
\draw [line width=2pt] (9,9)-- (14,9);
\draw [line width=2pt] (9,8)-- (14,8);
\draw [line width=2pt] (14,7)-- (9,7);
\draw [line width=2pt] (9,6)-- (12,6);
\draw [line width=2pt] (12,6)-- (12,9);
\draw [line width=2pt] (9,5)-- (11,5);
\draw [line width=2pt] (11,4)-- (9,4);
\draw [line width=2pt] (8,1)-- (8,9);
\draw [line width=2pt] (7,9)-- (8,9);
\draw [line width=2pt] (7,9)-- (7,1);
\draw [line width=2pt] (7,1)-- (8,1);
\draw [line width=2pt] (6,3)-- (6,9);
\draw [line width=2pt] (5,9)-- (5,3);
\draw [line width=2pt] (4,4)-- (4,9);
\draw [line width=2pt] (3,9)-- (3,6);
\draw [line width=2pt] (2,7)-- (2,9);
\draw [line width=2pt] (1,9)-- (1,7);
\draw [line width=2pt] (9,9)-- (9,3);
\draw [line width=2pt] (7,2)-- (8,2);
\draw [line width=2pt] (8,3)-- (7,3);
\draw [line width=2pt] (7,4)-- (8,4);
\draw [line width=2pt] (8,5)-- (7,5);
\draw [line width=2pt] (7,6)-- (8,6);
\draw [line width=2pt] (8,7)-- (7,7);
\draw [line width=2pt] (7,8)-- (8,8);
\draw [line width=2pt] (10,9)-- (10,3);
\draw [line width=2pt] (11,4)-- (11,9);
\draw [line width=2pt] (13,7)-- (13,9);
\draw [line width=2pt] (14,9)-- (14,7);
\draw [line width=2pt,color=ccqqqq] (7.5,1.5)-- (7.5,8.5);
\draw [line width=2pt,dash pattern=on 2pt off 4pt,color=ccqqqq] (7.5,8.5)-- (9.5,8.5);
\draw [line width=2pt,color=ccqqqq] (9.5,8.5)-- (13.5,8.5);
\draw [line width=2pt,color=ccqqqq] (9.5,3.5)-- (9.5,7.5);
\draw [line width=2pt,color=ccqqqq] (9.5,7.5)-- (13.5,7.5);
\draw [line width=2pt,color=ccqqqq] (10.5,4.5)-- (10.5,6.5);
\draw [line width=2pt,color=ccqqqq] (10.5,6.5)-- (11.5,6.5);
\draw [line width=2pt,color=qqwuqq] (5.5,3.5)-- (5.5,8.5);
\draw [line width=2pt,color=qqwuqq] (5.5,8.5)-- (1.5,8.5);
\draw [line width=2pt,color=qqwuqq] (4.5,4.5)-- (4.5,7.5);
\draw [line width=2pt,color=qqwuqq] (4.5,7.5)-- (1.5,7.5);
\begin{scriptsize}
\draw [fill=ccqqqq] (7.5,1.5) circle (2.5pt);
\draw [fill=ccqqqq] (7.5,8.5) circle (2.5pt);
\draw [fill=ccqqqq] (13.5,8.5) circle (2.5pt);
\draw [fill=ccqqqq] (9.5,8.5) circle (2.5pt);
\draw [fill=ccqqqq] (9.5,3.5) circle (2.5pt);
\draw [fill=ccqqqq] (9.5,7.5) circle (2.5pt);
\draw [fill=ccqqqq] (13.5,7.5) circle (2.5pt);
\draw [fill=ccqqqq] (10.5,4.5) circle (2.5pt);
\draw [fill=ccqqqq] (10.5,6.5) circle (2.5pt);
\draw [fill=ccqqqq] (11.5,6.5) circle (2.5pt);
\draw [fill=qqwuqq] (5.5,3.5) circle (2.5pt);
\draw [fill=qqwuqq] (5.5,8.5) circle (2.5pt);
\draw [fill=qqwuqq] (1.5,8.5) circle (2.5pt);
\draw [fill=qqwuqq] (1.5,7.5) circle (2.5pt);
\draw [fill=qqwuqq] (4.5,7.5) circle (2.5pt);
\draw [fill=qqwuqq] (4.5,4.5) circle (2.5pt);
\draw [fill=qqwuqq] (3.5,6.5) circle (2.5pt);
\draw (3.5,9.5) node[anchor=center] {\Large \reflectbox{$\pi$}};
\draw (11.5,9.5) node[anchor=center] {\Large $\pi$};
\end{scriptsize}
\end{tikzpicture}
\caption{Example the combinatorial proof of \eqref{eq:RefOdds1Psi} when $N=14$, $j=8$, and $\pi=(5,5,3,2,2,1)$. Then the partition $\pi_o$ is $(11,11,7,5,5,3,1,1)$ (read row by row). The partition we get after Sylvester bijection is $\pi_d = (13,10,9,7,4,1)$. The red lines (right-bending lines) are the odd-indexed parts of $\pi_d$ and the green lines (left-bending lines) are the even-indexed parts of $\pi_d$. Only the green lines are counted under the statistic $\E(\pi_d)$.}\label{figure_ferrers}
\end{figure}
\end{center}

Not only that Theorem~\ref{thm:RefinedOdds1} itself is elegant, it also has really elementary combinatorial corollaries as $N$ tends to infinity. 

\begin{corollary}\label{cor:Limit}
\begin{align}
\label{eq:RefOdds1PsiCor}\sum_{\substack{\pi\in\D\\\gamma(\pi) = j}} q^{\E(\pi)} &=  \frac{1}{(q;q)_j},\\
\label{eq:RefOdds1PhiCor}\sum_{\substack{\pi\in\P\\\gamma(\pi) = j}}q^{\E(\pi)} &= \frac{1}{(q;q)_j(q;q)_\infty}.
\end{align}
\end{corollary}

A combinatorial interpretation of \eqref{eq:RefOdds1PsiCor} is Theorem~\ref{cor:Intro}, and in the same spirit a combinatorial interpretation of \eqref{eq:RefOdds1PhiCor} is the theorem below.
\begin{theorem} The number of partitions $\pi$ where $\E(\pi)=n$ and $\gamma(\pi)=j$ is equal to the number of partitions of $n$ where parts are in two colors red and green such that red parts $\leq j$.
\end{theorem}

\section{Identities for weighted Rogers--Ramanujan partitions}\label{Sec6}

In \cite{Uncu}, the second author demonstrated how some of these theorems can be turned into weighted identities using multiplicative weights as well. In that paper it was shown that Schmidt's problem is equivalent to Alladi's weighted Rogers-Ramanujan count (see, \cite[(5.2)]{Uncu}).

\begin{theorem}\label{Uncu1} Let $\#(\pi)$ be the number of parts in $\pi$ and let $ \mathcal{R}\mathcal{R}$ be the set of partitions with gaps between parts $\geq 2$, then \begin{equation} \sum_{\pi\in\D}q^{\O(\pi)} = \sum_{\pi\in \mathcal{R}\mathcal{R}} \omega(\pi) q^{|\pi|}, \end{equation}where
\begin{equation}\label{eq:omega}\omega(\pi) := \lambda_{\#(\pi)}\cdot\prod_{i=1}^{\#(\pi)-1} (\lambda_{i}-\lambda_{i+1}-1).\end{equation}
\end{theorem}

Recently the same observation was made by Alladi \cite{AlladiSchmidt}.

It is clear that \eqref{eq:New_PsiN_rep} translates to the refinement of Theorem~\ref{Uncu1} with an extra bound on the largest part of partitions.

\begin{theorem} \label{New_weighted_sum} Let $ \mathcal{R}\mathcal{R}_{\leq N}$ be the set of partitions into parts $\leq N$ with gaps between parts $\geq 2$, then \begin{equation}\sum_{\pi\in\D_{\leq N}}q^{\O(\pi)} = \sum_{\pi\in \mathcal{R}\mathcal{R}_{ \leq N}} \omega(\pi) q^{|\pi|}.\end{equation}
\end{theorem}

Similarly, we can also interpret \eqref{eq:Odds1Psi} with the slight modification of the weight \eqref{eq:omega}.

\begin{theorem}\label{Uncu2} \begin{equation} \sum_{\pi\in\D_{\leq N}}q^{\E(\pi)} = \sum_{\pi\in \mathcal{R}\mathcal{R}_{\leq N-1}}\omega(\pi) q^{|\pi|} ,\end{equation}where
\begin{equation}\label{eq:omegaHat}\hat{\omega}(\pi) := (N-\lambda_1)\cdot\lambda_{\#(\pi)}\cdot\prod_{i=1}^{\#(\pi)-1} (\lambda_{i}-\lambda_{i+1}-1)\ \ \text{and}\ \ \hat{\omega}(\emptyset)=N+1.\end{equation}
\end{theorem}

There are $N+1$ partitions $\pi$ in $D_{\leq N}$ ($\emptyset$, $(k)$ for $1\leq k\leq N$) that yields $\E(\pi)=0$. This is the reasoning behind the special definition of $\hat{\omega}(\emptyset)$ in \eqref{eq:omegaHat}.

An example with $N=4$ is given for the above theorem. The sum of all the second columns in the first three blocks is equal to the sum of the weights in the last column.

\[
\begin{array}{cc|cc|cc||cc}
\pi\in\D_{\leq4} &  q^{\E(\pi)} & \pi\in\D_{\leq4} &  q^{\E(\pi)} & \pi\in\D_{\leq4} &  q^{\E(\pi)} &\pi\in\mathcal{R}\mathcal{R}_{\leq 3} &  \hat{\omega}(\pi)q^{|\pi|}\\
\hline
 (4)	& 1		& (4,2,1)	& q^2	& (2)	& 1		&	(3,1)	&	q^4		\\
 (4,3)	& q^3	& (4,3,2,1)	& q^4	& (2,1)	& q 	&	(3)	&	3q^3\\
 (4,2)	& q^2	& (3)		& 1	& (1)	& 1			&	(2)	&	4q^2	\\
 (4,1)	& q		& (3,2)		& q^2	& \emptyset	& 1	&	(1)	&	3q	\\
(4,3,2) & q^3	& (3,1)		& q	& 	& 				&	\emptyset	& 5	\\
(4,3,1)	& q^3	& (3,2,1)	& q^2	& 	& 			&
\end{array}
\]

\section{Acknowledgment} 
%%%%%%%%%%%%%%%%%%%%%%%%%%%%%%%%%%%%%%%%%%%%%%%%%%%%%%%%%%%%%%%%%%%%%%%%%%%%%%%%%%%%%%%%%%%%%%%%%%%%%

We would like to thank George E. Andrews for his kind interest and encouragement.

The second author thanks EPSRC and FWF for partially supporting his research through grants EP/T015713/1 and P-34501N, respectively.


\begin{thebibliography}{99}

\bibitem{AlladiWeighted} K. Alladi, \textit{Partition identities involving gaps and weights}, Trans. Amer. Math. Soc. \textbf{349} (1997), no. 12, 5001-5019.

\bibitem{AlladiSchmidt} K. Alladi, \textit{Schmidt-type theorems via weighted partition identities}, preprint.
   
\bibitem{theoryofpartitions}G. E. Andrews, \textit{The theory of partitions}, Cambridge Mathematical Library, Cambridge University Press, Cambridge, 1998. Reprint of the 1976 original. MR1634067 (99c:11126).

\bibitem{AndrewsKeith} G. E. Andrews, and W. J. Keith, \textit{Schmidt-type theorems for partitions with uncounted parts}, preprint \href{https://arxiv.org/abs/2203.05202}{arXiv:2203.05202}.

\bibitem{AndrewsPaule}  G. E. Andrews, and P. Paule, \textit{MacMahon’s partition analysis XIII: Schmidt type partitions and modular forms}, J. Number Theory, (2021), \href{https://doi.org/10.1016/j.jnt.2021.09.008}{https://doi.org/10.1016/j.jnt.2021.09.008}.

\bibitem{BU1} A. Berkovich, and A. K. Uncu, \textit{A new Companion to Capparelli's Identities}, Adv. in Appl. Math. \textbf{71} (2015), 125-137.

\bibitem{BU2} A. Berkovich, and A. K. Uncu, \textit{On partitions with fixed number of
even-indexed and odd-indexed odd parts}, J. Number Theory \textbf{167} (2016), 7-30.

\bibitem{WarnaarBerkovich} A. Berkovich, and S. O. Warnaar, \textit{Positivity preserving transformations for q-binomial coefficients}, Trans. Amer. Math. Soc. \textbf{357} (2005), no. 6, 2291-2351.

\bibitem{Bressoud} D. M. Bressoud, \textit{Proofs and Confirmations: The Story of the Alternating-Sign Matrix Conjecture}. Cambridge University Press, 1999.

\bibitem{BridgesUncu} W. Bridges, and A. K. Uncu, \textit{Weighted cylindric partitions}, preprint \href{https://arxiv.org/abs/2201.03047}{arXiv:2201.03047}.

\bibitem{Boulet} C. E. Boulet, \textit{A four-parameter partition identity}, Ramanujan J. \textbf{12} (2006), no. 3, 315-320.    

\bibitem{ChernYee} S. Chern and A. J. Yee, \textit{Diagonal Hooks and a Schmidt-Type Partition Identity}, Elect. J, of Comb. \textbf{29} vol 2, (2022), $\#$P2.10.

\bibitem{Fu} S. Fu, and J. Zeng \textit{A unifying combinatorial approach to refined little Göllnitz and Capparelli's companion identities}, Adv. in Appl. Math. \textbf{98} (2018), 127-154.

%\bibitem{GasperRahman} G. Gasper, and M. Rahman, \textit{Basic Hypergeometric Series}, Vol. 96. Cambridge university press, 2004.

\bibitem{Masao}  M. Ishikawa, and J. Zeng, \textit{The Andrews-Stanley partition function and Al-Salam-Chihara polynomials}, Discrete Math. \textbf{309} (2009), no. 1, 151-175.

\bibitem{Ji} K. Q. Ji. \textit{A Combinatorial Proof of a Schmidt Type Theorem of Andrews and Paule}, Elect. J, of Comb. \textbf{29} vol 1, (2022), $\#$P1.24.

\bibitem{LiYee} R. Li, and A.J. Yee, \textit{Schmidt type partitions}, preprint \href{https://arxiv.org/abs/2204.02535}{arXiv:2204.02535}.

\bibitem{Mork} P. Mork, \textit{Interrupted Partitions - Solution to Problem 10629} The Amer. Math. Monthly, Vol. \textbf{107}, No. 1 (Jan., 2000), pp. 87-87.

\bibitem{Schmidt} F. Schmidt, \textit{Interrupted Partitions - Problem 10629} The Amer. Math. Monthly, Vol. \textbf{104}, (1999), pp. 87-88.

\bibitem{Uncu} A K. Uncu, \textit{Weighted Rogers-Ramanujan Partitions and Dyson Crank},  Ramanujan J. \textbf{46}  (2018),  no. 2, 579-591.


   \end{thebibliography}
\end{document}